\documentclass[prodmode,acmtecs]{acmsmall} 

\usepackage[ruled]{algorithm2e}
\usepackage{amssymb,amsmath,amsfonts}
\usepackage{verbatim,graphics,url,comment,color,epstopdf}

\renewcommand{\algorithmcfname}{ALGORITHM}
\SetAlFnt{\small}
\SetAlCapFnt{\small}
\SetAlCapNameFnt{\small}
\SetAlCapHSkip{0pt}
\IncMargin{-\parindent}

\acmVolume{9}
\acmNumber{4}
\acmArticle{39}
\acmYear{2010}
\acmMonth{3}

\newcommand{\1}[1]{\text{ones}(#1)}
\newcommand{\0}[1]{\text{zeros}(#1)}
\newcommand{\rank}{\text{rank}}

\newtheorem{bound}{Bound}
\newcommand{\newalgname}[1]{%
  \renewcommand{\ALG@name}{#1}%
}

\newcounter{sAlgCount}
\newenvironment{sAlg}{
\begin{algorithm}[!htp]
\renewcommand{\algorithmcfname}{Sampling Method}
\stepcounter{sAlgCount}
\renewcommand{\thealgocf}{\arabic{sAlgCount}}
}{\end{algorithm}}

\newcounter{liAlgCount}
\newenvironment{liAlg}{
\begin{algorithm}[!htp]
\renewcommand{\algorithmcfname}{Leverage Score Distribution}
\stepcounter{liAlgCount}
\renewcommand{\thealgocf}{\arabic{liAlgCount}}
}{\end{algorithm}}

\newcounter{mtxGenAlg}
\newenvironment{mtxGenAlg}{
\begin{algorithm}[!htp]
\renewcommand{\algorithmcfname}{Matrix Generation Algorithm}
\stepcounter{mtxGenAlg}
\renewcommand{\thealgocf}{\arabic{mtxGenAlg}}
}{\end{algorithm}}

\newcounter{regularAlg}

\definecolor{gray}{RGB}{150,150,150}

\begin{document}

\markboth{T. Wentworth et al.}{kappa\_SQ, A \texttt{Matlab} package for randomized sampling of matrices with orthonormal columns.}

\title{kappa\_SQ: A \texttt{Matlab} package for randomized sampling of matrices with orthonormal columns}
\author{Thomas Wentworth
\affil{North Carolina State University}
Ilse Ipsen
\affil{North Carolina State University }}

\begin{abstract}
The kappa\_SQ software package is designed to assist researchers working on randomized row sampling.  The package contains a collection of \texttt{Matlab} functions along with a GUI that ties them all together
and provides a platform for the user to perform experiments.

In particular, kappa\_SQ is designed to do experiments related to the two-norm condition number of a sampled matrix, $\kappa(SQ)$, where $S$ is a row sampling matrix and $Q$ is a tall and skinny matrix with orthonormal columns.
Via a simple GUI, kappa\_SQ can generate test matrices, perform various types of row sampling, measure $\kappa(SQ)$, calculate bounds and produce high quality plots of the results.  All of the important codes are 
written in separate \texttt{Matlab} function files in a standard format which makes it easy for a user to either use the codes by themselves or incorporate their own codes into the kappa\_SQ package.\\

\end{abstract}

{}{}



\acmformat{Thomas Wentworth, and Ilse Ipsen, 2013. kappa\_SQ: A \texttt{Matlab} package for randomized sampling of matrices with orthonormal columns.}

\begin{bottomstuff}
This work is supported by NSF CISE CCF NCSU and GRANT NUMBERS.\\
Authors' address: Mathematics Department,
2311 Stinson Drive,
Box 8205, NC State University,
Raleigh, NC 27695-8205\\
Authors' e-mails: Thomas Wentworth: \texttt{thomas\_wentworth@ncsu.edu}, Ilse Ipsen: \texttt{ipsen@ncsu.edu}
\end{bottomstuff}

\maketitle
%
\section{Introduction}
We wrote the kappa\_SQ software package to assist us with our research on various algorithms for uniform row sampling \cite{ourPaper}.
In our research, a $m \times n$ matrix $Q$ with orthonormal columns and $m\geq n$ is sampled by a row sampling matrix $S$ to create the $c \times n$ sampled matrix $SQ$.
We then address the question, given $\eta>0$, what is the probability that rank$(SQ)=n$ and the two-norm condition number $\kappa(SQ) = \|SQ\|_2\|(SQ)^{\dagger}\|_2 \leq 1+\eta$?
This question is important due to its applications to randomized least squares solvers such as LSRN \cite{MSM11} and, in particular, the \textsl{Blendenpik} algorithm \cite{AMTol10}.\\
The \textsl{Blendenpik} algorithm uses randomized sampling to solve an overdetermined least-squares problem $\min_x\|Ax-b\|_2$ faster than LAPACK.  It starts by finding the $QR$ factorization, $Q_sR_s = SA$, of the randomly
sampled matrix $SA$ and then, if $SA$ has full colun rank, solves the preconditioned least squares problem
$\min_z\|AR_s^{-1}z - b\|_2$ via LSQR.
The solution to the original least squares problem can be found by solving a much smaller linear system with coefficient matrix $R_s$.
The key to this method is that if $\kappa(AR_s^{-1}) \approx 1$, then LSQR will converge quickly.\\

The connection between our work, kappa\_SQ and the \textsl{Blendenpik} algorithm is that if $SA$ is full rank, then $\kappa(SQ) = \kappa(AR_s^{-1})$.
This means that sampling rows from $A$ is, conceptually, the same as sampling rows from $Q$ and that $\kappa(AR_s^{-1})$ depends only
on the columns space of $A$ (and the sampling matrix).
Thus, it suffices to examine the behavior of $\kappa(SQ)$.\\

This code examines $\kappa(SQ)$ in two main ways.  First, it can perform numerical experiments where $\kappa(SQ)$ is measured.  And second, our code can plot bounds
for $\kappa(SQ)$.  In the literature, these bounds are often expressed in terms of two matrix properties that have been shown to be related to row sampling, 
leverage scores, and coherence.\\

Leverage scores were first introduced in 1978 by Hoaglin and Welsch \cite{HoagW78} to detect outliers when computing regression diagnostics.  They give a measurement of the distribution of the elements
in an orthonormal basis. 
The leverage scores of a matrix $A$ are defined in terms of any orthonormal basis, $Q$, for the column space of $A$.
\begin{definition}
The leverage scores of the real $m \times n$ matrix $A$ with $m \geq n$ are
$$\ell_j(A) = \ell_j(Q)\equiv \|e_j^TQ\|_2^2, \qquad 1\leq j\leq m.$$
\end{definition}
Since leverage scores are simply row norms from matrices with orthonormal columns, the inequality $0 \leq \ell_i(A) \leq 1$ holds and $\sum_{i=1}^{m}\ell_i(A) = n$.  If $\ell_i(A)=1$ then the $i'$th row
contains all of the information for a particular column.  On the other hand, if $\ell_i(A)=0$, then the $i'$th row of $A$ is zero and contains no data.
Thus, leverage scores give a quantification of the importance of each row with respect to sampling.
We use leverage scores as an input to both generate test matrices and to bound the condition number of a sampled matrix.
Our code for computing leverage scores is \texttt{leverageScores.m}.\\

In our work, coherence is simply the largest leverage score.

\begin{definition}[Definition 3.1 in \cite{AMTol10}, Definition 1.2
in \cite{CanR09}]\label{d_co}
The coherence of $A$ is 
\begin{eqnarray*}\label{e_coherence}
\mu(A) \equiv \max_{1\leq j\leq m}\ell_j(Q) = \max_{1\leq j\leq m}\|e_j^TQ\|_2^2.
\end{eqnarray*}
\end{definition}
Although coherence contains far less
information about a matrix than the leverage scores, it can still be useful in bounding the condition
number of a sampled matrix (see Bound \ref{t_cohbound})
and may be easier to estimate than leverage scores
Due to the properties of leverage scores, the inequality $n/m \leq \mu(A) \leq 1$ holds, and if $\mu(A) \approx n/m$, then $\ell_i(A) \approx n/m$.
Our code for computing coherence is \texttt{coherence.m}.\\

\section{kappa\_SQ Design}
Kappa\_SQ was designed to perform all of the computations from \cite{ourPaper} and output paper-ready plots.  It can assist researchers in the following ways.\\

First, the GUI for kappa\_SQ has been designed to assist the user set-up, perform and plot experiments on $\kappa(SQ)$.  There are two types of experiments,
the computation of (possibly probabilistic) bounds on $\kappa(SQ)$ and the computation of $\kappa(SQ)$ for a given or generated test matrix $Q$.
The GUI has also been coded to allow a user to easily incorporate their own codes by simply placing a properly formatted \texttt{Matlab} function in the ``boundsAndAlgorithms'' directory.\\

Second, kappa\_SQ includes a collection of codes for various algorithms and bounds pertaining to the field of randomized row sampling.  These codes
are all written as \texttt{Matlab} function files that can be used on their own or with the kappa\_SQ GUI.  The codes include functions for row sampling, test matrix generation, leverage score
distribution generation and functions to compute bounds for $\kappa(SQ)$.  The included codes are outlined in section \ref{algIntro}.\\

The kappa\_SQ codes can be broken up into two main groups, the GUI and ``Algorithm Codes.''
Below, we describe these codes and their functions.
\subsection{kappa\_SQ GUI}\label{guiSection}
The kappa\_SQ GUI is designed to produce plots of both numerical experiments, where $\kappa(SQ)$ is actually measured, and of bounds on $\kappa(SQ)$.
To perform a numerical experiment, kappa\_SQ will perform sampling on a matrix and then measure the condition number of the sampled matrix, $\kappa(SQ)$.  Occasionally,
the sampled matrix, $SQ$, is not full rank.  This event is termed a ``failure'' event and kappa\_SQ also keeps track of these.\\

\begin{figure}[!htp]
\centerline{\includegraphics[width=4in]{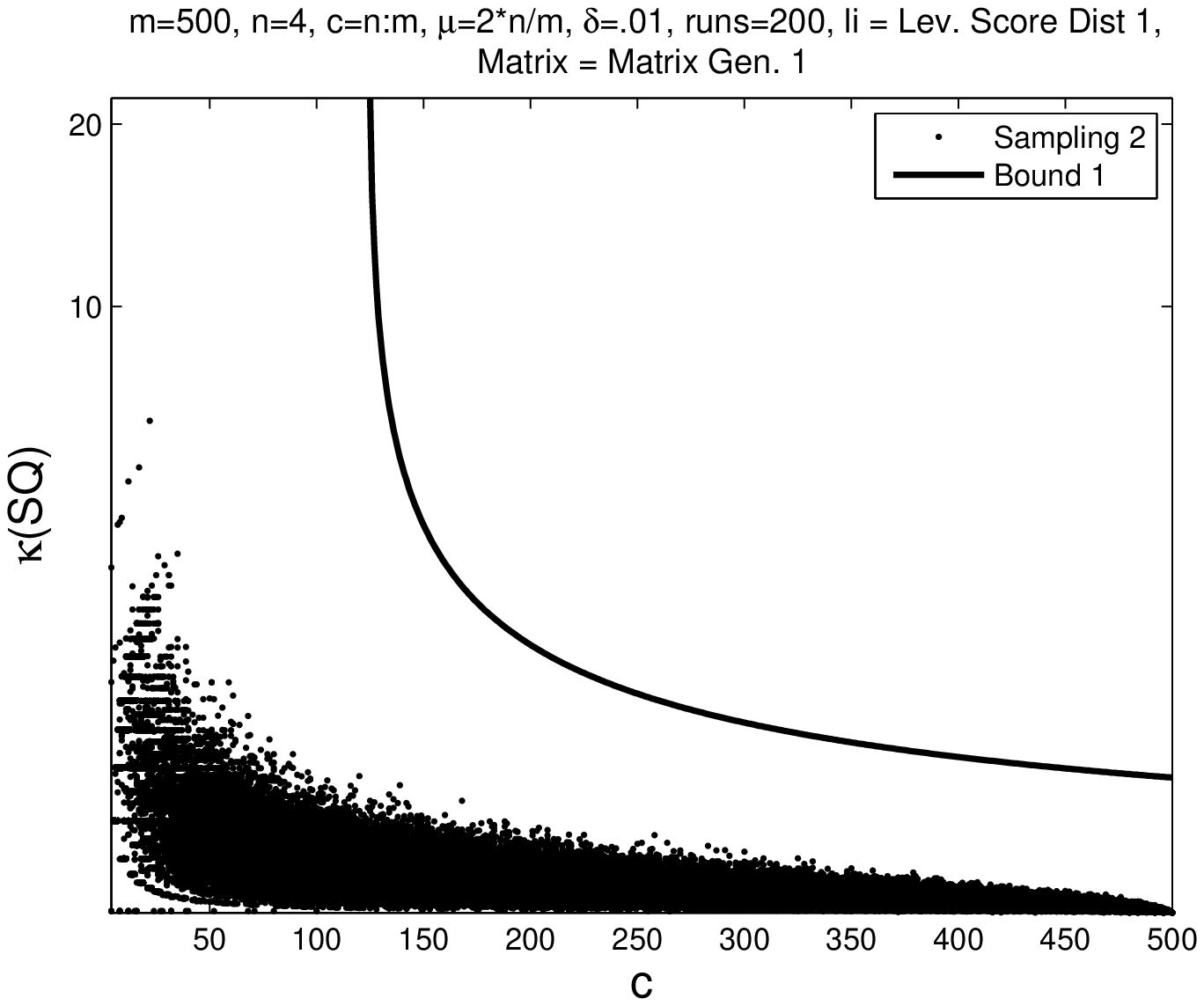}}
\caption{In this plot we show the results of a numerical experiment (triangles) and a bound on kappa\_SQ (line) that holds
with probability $1-\delta$.}
\label{exGui1}
\end{figure}

\begin{figure}[!htp]
\centerline{\includegraphics[width=4in]{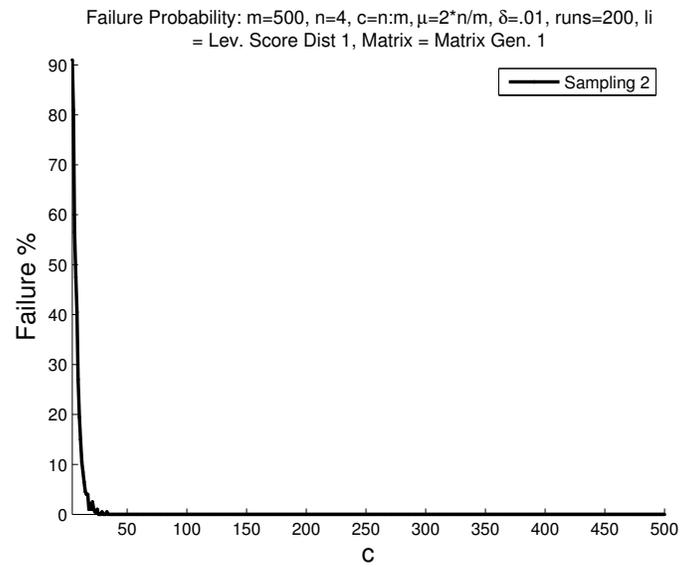}}
\caption{In this plot we show the failure rate of a numerical experiment on kappa\_SQ.}
\label{exGui2}
\end{figure}

When kappa\_SQ has completed its computations, it will output plots like those shown in Figures \ref{exGui1} and \ref{exGui2}.  For the moment, do not worry about
the specifics of each plot other than the following. The triangles in figure \ref{exGui1} show the measured condition number of the sampled matrix, $\kappa(SQ)$,
and the line plots a bound on $\kappa(SQ)$.  For the sections of the domain where a line is not plotted, the bound does not apply.
All of the bounds included with $\kappa(SQ)$ are probabilistic bounds and therefore only hold with probability at least $1-\delta$.
Therefore, at least $100(1-\delta)\%$ of the measured $\kappa(SQ)$ should be below the line.  In Figure \ref{exGui2}, the ``failure rate,'' the percent of numerical experiments
that resulted in a failure event, is plotted.  Despite the fact that for many experiments the failure rate will be $0\%$ for most values, it is still an important quantity because 
figure \ref{exGui1} only plots the the ``good'' events (where $SQ$ has full column rank).\\


\begin{figure}[!htp]
\centerline{\includegraphics[width=4in]{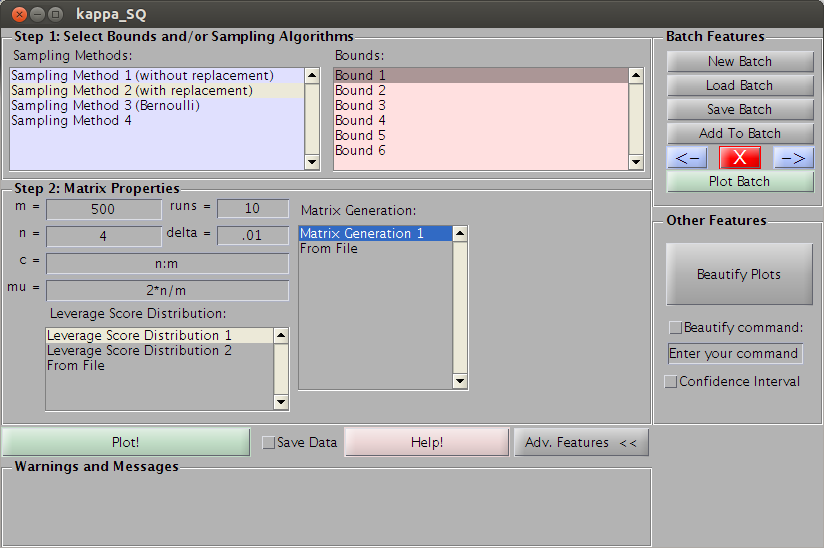}}
\caption{kappaSQ GUI with advanced features shown.}
\label{gui1}
\end{figure}

When a user executes \texttt{kappaSQ.m}, he or she is presented with the kappa\_SQ GUI (see Figure \ref{gui1}).
It is broken up into sections; the first two
correspond to inputs the user must provide to produce a plot and the last two are for plot editing and
and to display important information. Below we describe these sections and, in the process, how to use
the GUI to produce plots.\\

\textbf{``Step 1: Select Bounds and/or Numerical Experiments.''} This section allows the user to select
what he or she would like to plot.  The first listbox contains possible bounds that the user can plot.
The second listbox contains various sampling methods.  In kappa\_SQ, numerical experiments are first defined in terms of what sampling
algorithm the user would like to use.  By selecting a sampling method, the user is telling kappa\_SQ that he or she wishes to do a numerical experiment
with that sampling method.  Selecting multiple sampling methods will perform multiple experiments.\\

\textbf{``Step 2: Matrix Properties / Parameters.''}  This section allows the user to provide the required inputs.
Only inputs that are required will be visible.  As an example, if the user chose to plot Bound \ref{t_cohbound}, then this section would ask for the user
to provide values for $m, n, \mu, \delta$ and $c$ as those values are required to compute Bound \ref{t_cohbound}.
In addition, either $c$ or $\mu$ must be a vector and whichever is a vector will be placed on the $x$-axis.\\

Of particular interest in this section are the ``Matrix Generation'' and ``li'' ($\ell_i$) inputs. When a test matrix is required (ex: when running a numerical experiment), kappa\_SQ will generate a matrix
using the algorithm specified in this listbox.  Similarly, when a leverage score distribution is required, kappa\_SQ will generate one by the method specified in the ``li'' listbox.\\

\textbf{``Plot Button.''} Once the user has completed steps 1 and 2, he or she may click the plot button to run and plot the experiment.\\

\textbf{``Help!''} This button will open the kappa\_SQ help file.  This file includes a FAQ section and a list of included functions.\\

\textbf{Adv. Features: ``Batch Features.''} This section can be viewed by clicking on the ``Adv. Features'' button.  After performing steps 1 and 2, the user can instead
add the current experiment to a batch of jobs to be run later in serial.  This is particularly useful if the chosen experiments require a long time to run, or if the user has many experiments to run.
In addition, if the user keeps all of their experiments defined in a batch file, they can easily repeat all of the experiments for their work.\\

\begin{figure}[!htp]
\centerline{\includegraphics[width=4in]{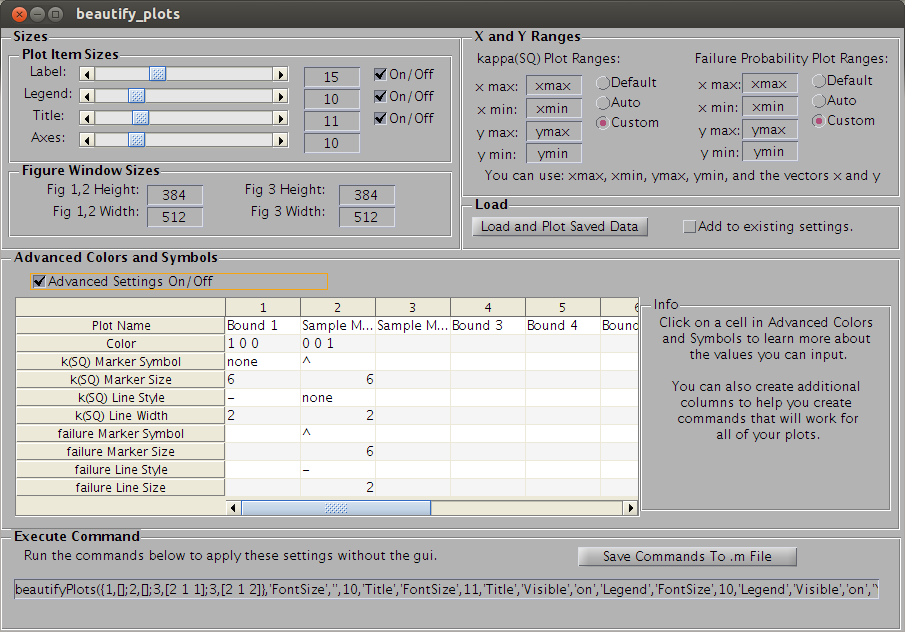}}
\caption{Beautify Plots GUI.}
\label{gui2}
\end{figure}

\textbf{Adv. Features: ``Other Features.''}  This section contains a button called ``Beautify Plots'' which will open the plot editing window shown in Figure \ref{gui2}.
This window assists the user with modifying many of the common plot settings and creating a script that will apply these settings to future plots.  In addition, the plot
editing window will generate a command which will apply these settings without the GUI.  KappaSQ can be set to run this command for all future plots by checking
the ``beautify command'' checckbox and entering the command.  Thus, the user only needs to set up their plots once.  Finally, this section also has an option to plot
a standard confidence interval for the failure probability.


Below we will discuss the various included codes that the kappa\_SQ GUI uses in the above sections to produce plots.

\subsection{Algorithm Codes}\label{algIntro}
In this section, we describe the various algorithms and bounds that are included in the kappa\_SQ package.
These algorithms can be broken up into four main groups, bounds, matrix generation, sampling methods and leverage score distributions.

\subsubsection{Sampling methods}  \label{sampMethods}
We include four different row sampling methods from the literature, Sampling without Replacement, Sampling with Replacement, Bernoulli Sampling, and Sampling Proportional to Leverage Scores.  In our recent paper \cite{ourPaper}, the first three of these sampling methods are described by constructing a sampling matrix $S$ such that $SQ$ is the sampled matrix.
In the kappa\_SQ package, we instead code these algorithms to compute $B=\equiv SQ$ directly.
Each sampling method inputs the initial matrix, $Q$, and the desired (or desired expected) number of rows to be sampled, $c$, and outputs the sampled matrix, $SQ$.\\

\textbf{Sampling Method \ref{alg:without}, Sampling Without Replacement.} This algorithm samples exactly the desired number of rows such that no row is sampled more than once by sampling
uniformly from the $m!/(m-c)!$ possible permutations of $c$ rows.  We implement this by first randomly permuting all of the rows
and then sampling the first $c$ rows.  In the algorithm below we use the term random permutations. A permutation $\pi_1,\ldots, \pi_m$ of the integers $1,\ldots, m$
is a \textit{random permutation}, if it is equally likely to be one of $m!$
possible permutations \cite[pages 41 and 48]{MitzUpf}.  The matlab command \texttt{randperm(m)} will generate a random permutation of the integers $1,\ldots, m$.  Our code for this sampling method is \texttt{Sample\_randperm.m}.\\

\begin{sAlg}
\SetAlgoNoLine
\KwIn{A $m\times n$ matrix $Q$ and an integer $c$, such that $1\leq c\leq m$.}
\KwOut{A $c\times m$ sampled matrix $B$.}
$v = \texttt{randperm}(m)$\;
$s = v(1:c)$\;
$B=\sqrt{m/c}\>Q(s,:)$\;
\caption{ \textsl{Sampling Without Replacement}, \protect\cite{GT11,GrN10} \label{samp:SWoR} }
\label{alg:without}
\end{sAlg}

\textbf{Sampling Method \ref{alg:with}, Sampling With Replacement (Exactly(c)).} This algorithm samples exactly the desired number of rows with a uniform probability distribution and with replacement.
Our code for implementing this sampling method is
\texttt{Sample\_exactlyC.m}.\\
\begin{sAlg}
\SetAlgoNoLine
\KwIn{A $m\times n$ matrix $Q$ and an integer $c$, such that $1\leq c\leq m$.}
\KwOut{A $c\times m$ sampled matrix $B$.}
Let $\pi_1,\ldots,\pi_c$ be integers uniformly sampled from$\left\{1, \ldots, m\right\}$ with replacement\;
$s = \begin{bmatrix}\pi_1, & \ldots, & \pi_c\end{bmatrix}$\;
$B=\sqrt{m/c}\>Q(s,:)$\;
\caption{\textsl{Sampling With Replacement, \protect\cite{DKM06,DMMS10}} \label{samp:SWR}
}
\label{alg:with}
\end{sAlg}

\textbf{Sampling Method \ref{alg:ber}, Bernoulli Sampling.} In this sampling method, each row is either sampled, with probability $c/m$, or not sampled with probability $1-c/m$.  Thus, whether or not each row is sampled
is an independent Bernoulli trial and the expected total number of rows sampled is $c$.  Our implementation of this algorithm differs slightly from \cite{ourPaper}.
In our paper, rows that are not sampled are set to $0$, while in our code they are removed.  Removing the zero rows is more memory efficient,
avoids unnecessary matrix-matrix multiplications and does not affect $\kappa(SQ)$.  Our code for this algorithm is
\texttt{Sample\_bernoulli.m}.\\

\begin{sAlg}
\SetAlgoNoLine
\KwIn{A $m\times n$ matrix $Q$ and an integer $c$, such that $1\leq c \leq m$.}
\KwOut{A $c\times m$ sampled matrix $B$.}
Let $\pi$ be a $m \times 1$ vector of $m$ independent realizations of a boolean random variable with success probability $c/m$\;
Let $s$ be a $\hat{c}\times 1$ vector containing the indices where $\pi_i=1$, where $1\leq i\leq m$ and $\hat{c} = $ the number of nonzero entries in $\pi$\;
$B=\sqrt{m/c}\>Q(s,:)$\;
\caption{\textsl{Bernoulli Sampling, \protect\cite{AMTol10,GT11,GrN10}} \label{samp:Ber}
}
\label{alg:ber}
\end{sAlg}

\textbf{Sampling Method \ref{alg:samp_li}, Sampling Proportional to Leverage Scores.}
In this sampling method, $c$ rows are sampled with probability $\ell_i(Q)/n$ with replacement.
Our code for implementing this sampling method is \texttt{Sample\_leverageScores.m}\\
\begin{sAlg}
\SetAlgoNoLine
\KwIn{A $m\times n$ matrix $Q$ and an integer $c$, such that $1\leq c\leq m$ and the leverage scores $\ell(Q)$.}
\KwOut{A $c\times m$ sampled matrix $B$.}
Let $\pi_1,\ldots,\pi_c$ be integers sampled from $\left\{1, \ldots, m\right\}$ with probabilities $\left\{\ell_1(Q)/n, \ldots, \ell_m(Q)/n\right\}$ and  replacement\;
$s = \begin{bmatrix}\pi_1, & \ldots, & \pi_c\end{bmatrix}$\;
$B=\sqrt{m/c}\>Q(s,:)$\;
\caption{\textsl{Sampling Proportional to Leverage Scores} \label{samp:SPLi}
}
\label{alg:samp_li}
\end{sAlg}

\subsubsection{Bounds}\label{bound1}
We include codes for the two probabilistic bounds for $\kappa(SQ)$ from our recent paper \cite{ourPaper} and four other weaker bounds that were included
in the first virsion of our paper \cite{ourOldPaper}.

\textbf{Bound \ref{t_cohbound}, Coherence based bound.} This bound is expressed in terms of coherence and comes from a matrix Chernoff concentration inequality \cite[Corollary 5.2]{tropp11}.
It applies to all of the sampling methods \ref{sampMethods} except for Sampling Proportional to Leverage Scores (Sampling Method \ref{alg:samp_li}). Our code for this bound is \texttt{Bound\_muBound.m}.\\
\begin{bound}[\protect{\cite[Theorem 4.1]{ourPaper}}]\label{t_cohbound}
Let $Q$ be a real $m\times n$ matrix with $Q^TQ=I_n$ and coherence~$\mu$.
Let $SQ$ be a sampling matrix produced by 
Algorithms~\ref{alg:without}, \ref{alg:with}, or \ref{alg:ber}
with $n\leq c\leq m$.
For $0<\epsilon<1$ and $f(x)\equiv e^x(1+x)^{-(1+x)}$ define
\begin{equation}\delta\equiv n\left(f(-\epsilon)^{c/(m\mu)}+
f(\epsilon)^{c/(m\mu)}\right).\label{cohBound:deltaEqn}\end{equation}
If $\delta<1$, then with probability at least $1-\delta$ we have
$\text{rank}(SQ)=n$ and
$$\kappa(SQ)\leq \sqrt{\frac{1+\epsilon}{1-\epsilon}}.$$
\end{bound}

\textbf{Bound \ref{t_levbound}, Leverage score based bound.} This bound is based on leverage scores and only applies to sampling without replacement (Sampling Method \ref{alg:without}).
It is based on a matrix Bernstein concentration inequality \cite[Theorem 4]{Recht11}\cite[Theorem 8.2]{ourPaper}. Our code for this bound is \texttt{Bound\_leverageScoresBound.m}.\\
\begin{bound}[\protect{\cite[Theorem 5.2]{ourPaper}}]\label{t_levbound}
Let $Q$ be a $m\times n$ real matrix with $Q^TQ=I_n$, leverage scores
$\ell_j(Q)$, $1\leq j\leq m$, and coherence $\mu$. Let $L$ be a diagonal matrix such that $L_{j,j} = \ell_j(Q)$. Let $S$ be a sampling
matrix produced by Algorithm~\ref{alg:with} with
$n\leq c\leq m$. For  $0<\epsilon<1$ set
$$\delta \equiv 2n\exp\left(-\tfrac{3}{2}\>
\frac{c \epsilon^2}{m\>(3\|Q^TLQ\|_2+\epsilon\mu)}\right).$$
If $\delta<1$, then with probability at least $1-\delta$ we have
$\text{rank}(SQ)=n$ and 
$$\kappa(SQ) \leq \sqrt{\frac{1+\epsilon}{1-\epsilon}}.$$
\end{bound}

\textbf{Bound \ref{t_b1}, Weaker coherence based bound.} This bound is based on a probabilistic two-norm bound for a
Monte Carlo matrix multiplication algorithm that samples according to Sampling Method \ref{alg:with} \cite[Theorem 4]{DMMS10}. Our code for this bound is \texttt{weakerBound\_1.m}.\\
\begin{bound}[\protect{\cite[Theorem 3.2]{ourOldPaper}}]\label{t_b1}
Given $0<\epsilon<1$ and $0<\delta<1$, let $Q$ be a $m\times n$ real matrix with $Q^TQ=I_n$ and coherence $\mu$.
Let $c$ be an integer so that
$$\min\left\{n,\zeta\ln{\left(\zeta/\sqrt{\delta}\right)}\right\}
\leq c\leq m,\qquad
\text{where} \qquad \zeta\equiv \frac{96 m\,\mu}{\epsilon^2}.$$
If $S$ is a $c\times m$ matrix produced by Sampling Method \ref{alg:with}
with uniform probabilities  $p_k=1/m$, $1\leq k\leq m$, then with probability 
at least $1-\delta$, we have $\rank(SQ)=\rank(M_s)=n$ and
$$\kappa(SQ)=\kappa(AR_s^{-1})\leq \sqrt{\frac{1+\epsilon}{1-\epsilon}}.$$
\end{bound}

\textbf{Bound \ref{c_b1a}, Weaker coherence based bound.}
This bound is based on a special case of the noncommutative
Bernstein inequality \cite[Theorem 4]{Recht11} and applies to sampling Sampling Method \ref{alg:with}. Our code for this bound is \texttt{weakerBound\_3.m}.\\
\begin{bound}[\protect{\cite[Corollary 3.10]{ourOldPaper}}]\label{c_b1a}
Given $c\geq n$ and $0<\delta<1$, let $Q$ be a $m\times n$ real matrix with $Q^TQ=I_n$ and coherence $\mu$.
Let $\rho\equiv\frac{2}{3}\ln(2n/\delta)$
and
$$\epsilon_1\equiv \frac{\mu m}{2c}\>
\left(\rho +\sqrt{\frac{12 c \rho}{m\mu} + \rho^2}\right).$$
Let $S$ be a $m \times m$ matrix produced by Algorithm~\ref{alg:with}.
If $\epsilon_1<1$ then with probability at least $1-\delta$, we have
$\rank(SQ)=n$ and 
$$\kappa(SQ) \leq \sqrt{\frac{1+\epsilon_1}{1-\epsilon_1}}.$$
\end{bound}

\textbf{Bound \ref{t_b2}, Weaker coherence based bound.} This bound is based on a Frobenius norm bound for 
a Monte Carlo matrix multiplication algorithm that samples according
to Sampling Method \ref{alg:with} and applies to sampling Sampling Method \ref{alg:with}. Our code for this bound is \texttt{weakerBound\_4.m}.\\
\begin{bound}[\protect{\cite[Theorem 3.5]{ourOldPaper}}]\label{t_b2}
Given $0<\delta<1$ and $c\geq n$, let $Q$ be a $m\times n$ real matrix with $Q^TQ=I_n$ and coherence $\mu$.
Let
$$\epsilon_2\equiv \sqrt{\frac{mn\>\mu}{c}} + 
m\>\mu\sqrt{\frac{8 \log(1/\delta)}{c}}.$$
Let $S$ be a $c\times m$ matrix produced by Algorithm~\ref{alg:with}
with uniform probabilities $p_k=1/m$, $1\leq k\leq m$. 
If $\epsilon_2<1$, then with probability at least $1-\delta$, 
we have $\rank(SQ)=n$ and
$$\kappa(SQ)\leq \sqrt{\frac{1+\epsilon_2}{1-\epsilon_2}}.$$
\end{bound}

\textbf{Bound \ref{c_b3aux}, Weaker coherence based bound.} This bound is again based on the noncommutative Bernstein
inequality in \cite[Theorem 4]{Recht11} and applies to sampling Sampling Method \ref{alg:ber}. Our code for this bound is \texttt{weakerBound\_6.m}.\\
\begin{bound}[\protect{\cite[Corollary 4.3]{ourOldPaper}}]\label{c_b3aux}
Given $m\geq n$, $0<\gamma<1$ and $0<\delta<1$, let $Q$ be a $m\times n$ real matrix with $Q^TQ=I_n$ and coherence $\mu$.
Let $\rho\equiv\frac{2}{3}\ln(2n/\delta)$
and
$$\hat{\epsilon}_3\equiv \frac{\mu}{2}
\left(\phi \rho +\sqrt{\frac{1-\gamma}{\gamma}\>12m \rho +\phi^2 \rho^2}\right),
\qquad 
\phi=\begin{cases} 1 & \text{if $\gamma\geq 1-\gamma$} \\

\frac{1-\gamma}{\gamma} & \text{if $1-\gamma >\gamma$}
\end{cases}$$
Let $S$ be a $m \times m$ matrix produced by Algorithm~\ref{alg:ber}.
If $\hat{\epsilon}_3<1$ then with probability at least $1-\delta$, we have
$\rank(SQ)=n$ and 
$$\kappa(SQ) \leq 
\sqrt{\frac{1+\hat{\epsilon}_3}{1-\hat{\epsilon}_3}}.$$
\end{bound}

\subsubsection{Leverage score distribution}\label{lidist1}


%
%

We include code for two functions which define leverage score distributions.\\

\textbf{Leverage Score Distribution \ref{alg:liOneBig}, Good leverage score distribution.} The first function is designed to be an ideal case for row sampling.  It outputs a leverage score distribution with one leverage
score is set equal to the coherence and the remaining leverage scores all identical.
Thus, most rows are equally ``important'' and uniform row sampling should work well.
The code for this algorithm is \texttt{liDist\_oneBig.m}.\\

\begin{liAlg}
\caption{\textsl{Good leverage score distribution \protect\cite[Algorithm 6.2]{ourPaper}}}
\SetAlgoNoLine
\KwIn{Integers $m$ and $n$ such that $m\geq n \geq 1$, and desired coherence $\mu$.}
\KwOut{A $m\times 1$ vector, $\ell$ of leverage scores such that $\max{\ell} = \mu$.}
$\ell=[\mu; \1{m-1,1}(n-mu)/(m-1) ]$\;
\label{alg:liOneBig}
\end{liAlg}

\textbf{Leverage Score Distribution \ref{alg:liManyBig}, Bad leverage score distribution.} The second function is designed to be a particularly bad case for row sampling. It outputs a leverage score distribution with the maximal number of entries set equal to the coherence and at most one additional non-zero entry.
Matrices with this leverage score distribution will have the maximal number of zero rows for the given coherence.  Zero rows are bad for uniform row sampling since zero rows contain no information.
The code for this algorithm is \texttt{liDist\_manyBig.m}.\\

\begin{liAlg}
\caption{\textsl{Bad leverage score distribution \protect\cite[Algorithm 6.3]{ourPaper}}}
\SetAlgoNoLine
\KwIn{Integers $m$ and $n$ such that $m\geq n \geq 1$, and desired coherence $\mu$.}
\KwOut{A $m\times 1$ vector, $\ell$ of leverage scores such that $\max{\ell} = \mu$.}
$\tilde{m} = \left\lfloor n/\mu \right\rfloor$\;
\eIf{$\tilde{m}<m$}{
$\ell=[\mu\1{\tilde{m},1}; n-\tilde{m}\mu, \0{m-\tilde{m},1} ]$\;
}
{
$\ell = \mu\1{\tilde{m},1}$
}

\label{alg:liManyBig}
\end{liAlg}

\subsubsection{Test matrix generation}\label{mtxgen1}
Included in kappa\_SQ, we provide code for a deterministic matrix generation algorithm.\\

\textbf{Matrix Generation \ref{alg:mtxGen}, Deterministic matrix generation algorithm.} This algorithm inputs the desired matrix dimensions and leverage scores
and outputs a test matrix, $Q$, with these properties.
To compute $Q$, the algorithm applies $m-1$ Givens rotations to the matrix $Q=\begin{bmatrix}I_n \qquad \0{n,m-n }\end{bmatrix}^T$.  Each Givens rotation alters the leverage scores of two rows such that
at least one of the two rows has the desired leverage score.  The reason that only $m-1$ Givens rotations are required is that the leverage scores sum to $n$ and thus the final leverage sore is determined by the other $m-1$
leverage scores.  We note here that this algorithm is a transposed version of \cite[Algorithm 3]{DHST05} and that the Givens rotations are computed from
numerically stable expressions \cite[section 3.1]{DHST05}. The code for this algorithm is \texttt{mtxGen\_li.m}.\\
\begin{mtxGenAlg}
\SetAlgoNoLine
\KwIn{Integers $m$ and $n$ such that $m\geq n \geq 1$, and a $m \times 1$ vector $l$ of the desired leverage scores.}
\KwOut{A $m\times n$ matrix $Q$ with orthonormal columns and the desired leverage scores.}
$Q=\begin{bmatrix}I_n & \0{n,m-n}\end{bmatrix}^T$\;
$[l,I]$=sort($l$);\% Sort and store original order\;
$i=m-n$\;
$j=m-n+1$\;
\For{$dummyVar=1:m-1$}{
\eIf{ $\left|l_i - \|e_i^TQ\|_2\right| < \left|l_j - \|e_j^TQ\|_2\right|$}
{ Rotate rows $i$ and $j$ of $Q$ so that $\|e_i^TQ\|_2^2=l_i$\;
$i=i-1$\;}
{ Rotate rows $i$ and $j$ of $Q$ so that $\|e_j^TQ\|_2^2=l_j$\;
$j=j+1$\;}
}
$Q(I,:)=Q;$\% Undo sorting\;
\caption{\textsl{Matrix Generation \ref{alg:mtxGen}} 
}
\label{alg:mtxGen}
\end{mtxGenAlg}

\subsection{Other Functions}
We also include two simple functions to assist with choosing nice, aesthetically pleasing, ranges for $c$ and $\mu$ named \texttt{logPoints.m} and \texttt{logPointsDouble.m}, respectivelly.
These functions produce ranges that are more heavily weighted towards the smaller end of the desired range.  Most of the interesting action in the final plots occurs near smaller $c$ or $\mu$ values
and, in addition, larger values of $c$ are more computationally expensive.
We describe these functions in the kappa\_SQ help file which can be accessed by pressing the help button in the gui (see Section \ref{guiSection}).


\section{Examples}
Here we show a few examples of some of the ways that the kappa\_SQ GUI can be used.\\

\textbf{Example 1:} In this example, we show how to perform a basic experiment with the GUI.  In this experiment we compare Bound \ref{t_cohbound} to a numerical experiment
with sampling Sampling Method \ref{alg:with}.  Since Bound \ref{t_cohbound} applies to this sampling method, the results should show that at least $100(1-\delta)\%$ of the measured $\kappa(SQ)$ are less than
the bound.  In order for kappa\_SQ to run a numerical experiment, it must have a test matrix to work on.  Here, we will chose to generate
a test matrix with Sampling Method \ref{alg:mtxGen} and leverage scores defined by Sampling Method \ref{alg:liOneBig}.\\

\begin{figure}[!htp]
\includegraphics[width=4in]{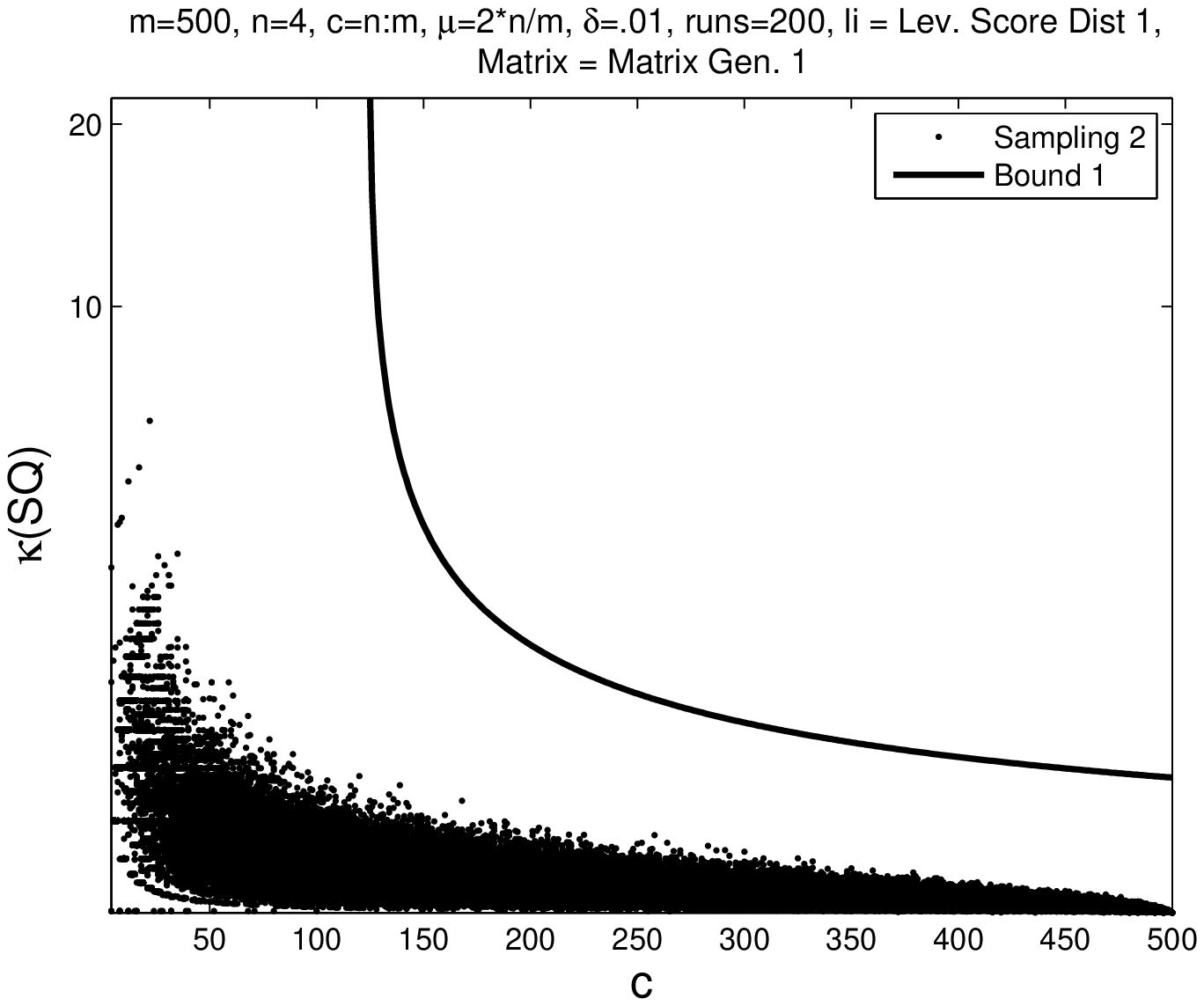}
\caption{Resulting plots produced after clicking on the plot-button shown in Figure \ref{ex2}.  The solid line shows Bound \ref{t_cohbound} and the triangles show the results of the numerical experiments
with sampling Sampling Method \ref{alg:with} and a matrix generated by Sampling Method \ref{alg:mtxGen}.

Here $Q$ is a matrix generated by
Algorithm~\ref{alg:mtxGen} with orthonormal 
columns, $m=10,000$, $n=4$, coherence $\mu=20n/m$, .
Left panel: Horizontal coordinate axes represent amounts of sampling 
$n\leq c\leq 10,000$. Vertical coordinate axes represent
condition numbers $\kappa(SQ)$; the maximum is 10.
Right panels: Horizontal coordinate axes represent amounts of sampling that
give rise to numerically rank deficient matrices $SQ$.
Vertical coordinate axes represent percentage of
numerically rank deficient matrices.
}
\label{ex2figs}
\end{figure}

To set up the experiment, start by selecting the desired bound and sampling method.  Then, move on to step 2 and input the following values, \texttt{m=500}, \texttt{n=4}
\texttt{c=n:m}, \texttt{mu=2n/m}, \texttt{runs=10}, and \texttt{delta=.01}.  For the ``Matrix Generation'' listbox, select ``Matrix Generation \ref{alg:mtxGen}'' (Sampling Method \ref{alg:mtxGen}).
This will cause the leverage score listbox to appear since this matrix generation algorithm requires a leverage score distribution.  Select ``Leverage Score Distribution 1'' (Sampling Method \ref{alg:liOneBig})
for the leverage score distribution.
In figure \ref{gui1} we show how the GUI should at this point.
When ready, click the plot button to begin the experiment. We show the resulting plot in figure \ref{ex2figs}.\\


\textbf{Example 2:} Here we show how to set up a kappa\_SQ batch to perform multiple experiments in serial.  To create a new batch, first click the ``Adv. Features'' button
to expand the GUI and then click on the ``New Batch'' button to start a new batch.
Next, set up an experiment by performing steps 1 and 2 as described in the first example.  Then, instead of clicking the plot button, click the ``Add to Batch'' button.
This will add the current experiment to the batch.  Repeat this process for the remaining experiments.
The user may use the arrow buttons to navigate and the ``X'' button to delete previously entered experiments.
When ready, click the ``Save Batch'' button to save the current experiments to a file and then the ``Run Batch" button to have kappa\_SQ run all of the experiments.
Plot images will be saved automatically with a file name based on the batch file name and their job number.\\

\textbf{Example 3:} Here we show how the built in plot editing tools can be used to expedite plot editing, how to create a script that will apply these settings to future plots and set up
kappa\_SQ to run that script after every experiment.  To start, run an experiment as described in example 1.  Then, click the ''Adv. Features`` button to expand the GUI
and then click on the ``Beautify Plots'' button. This will open the plot editing window. (See figure \ref{gui2}).\\

This window allows the user to easily edit many different plot settings.  While using this window, any changes will instantly be applied to the plots, so we suggest positioning the plots and the GUI window such
that they can all be seen. When done editing, press the button labeled ''Save Commands To .m File`` to create a \texttt{.m}
file that will apply these settings to future plots.  To have kappa\_SQ apply these plot settings automatically to all new plots, write the command for this file in the box labeled ''Enter your command here``
in the main GUI window and check the ''Beautify Command`` checkbox. (See figure \ref{gui1}).

\section{Conclusions}

The kappa\_SQ package is designed to assist researchers examine the behavior of $\kappa(SQ)$.
The package includes codes for generating matrices with specific leverage score distributions, generating a two specific leverage score distributions,
four types of row sampling methods and computing bounds on $\kappa(SQ)$.  These codes can be used on their own, or with the kappa\_SQ
GUI which is capable of setting up and running numerical experiments,
computing bounds, and producing quality plots with the help of custom plot-editing tools.  In addition, the GUI has been designed to detect properly formatted
\texttt{Matlab} function files which allows the user to incorporate their own codes into the GUI.
\bibliographystyle{ACM-Reference-Format-Journals}
\bibliography{biblio}

\section{Appendix}
\subsection{Notation}
We use the following notation through out this paper.
\begin{itemize}
 \item $m$, $c$ and $n$ are integers such that $m\geq c \geq n>0$.
 \item $\|\cdot\|_2$ denotes the standard $2$-norm.
 \item $A^T$ denotes the transpose of $A$.
 \item $e_i$ denotes the canonical vector with a 1 in the $i^\text{th}$ position and zeros everywhere else.
 \item $A$ is a $m \times n$ full column rank matrix.
 \item $Q$ is a $m\times n$ matrix with orthonormal columns that span the column space of $A$.
 \item $S$ is a $c \times n$ random row sampling matrix.
 \item $\kappa(A)\equiv \|A\|_2\|A^\dagger\|_2$ denotes the two-norm condition number of a $m \times n$ full column rank matrix $A$, where $A^\dagger$ is the Moore-Penrose inverse.
 \item $I_k = \begin{pmatrix} e_1 & \ldots & e_k \end{pmatrix}$ denotes the $k\times k$ identity matrix.
 \item $\1{m\times n}$ denotes the $m\times n$ matrix of all ones.
 \item $\0{m\times n}$ denotes the $m \times n$ matrix of all zeros.
 \item $\mu(A)$ denotes the coherence of $A$.
 \item $\ell(A)$ denotes the $m \times 1$ vector containing the leverage scores of $A$, and $\ell_i(A)$ denotes the $i^{\text{th}}$ leverage score of $A$.
 \item $L$ is a diagonal matrix with the leverage score of $Q$ on the diagonal.
 \item $\delta$ is a number such that $0<\delta  <1$ and is referred to as the failure probability.
 
\end{itemize}

\end{document}